\documentclass[a4paper, 10pt, reqno]{amsart}
\usepackage{latexsym}
\usepackage{fancyhdr}
\usepackage{amssymb}
\usepackage{array}
\usepackage{enumerate}
\usepackage{verbatim}
\usepackage{xspace}
\usepackage{exscale}
\usepackage[ansinew]{inputenc}
\usepackage{mathrsfs}
\usepackage[all]{xy}
\usepackage{tabularx}
\usepackage[final]{graphicx}

\theoremstyle{plain}
\newtheorem*{lemma*}{Lemma}

\newtheorem*{theorem*}{Theorem}

\newtheorem*{proposition*}{Proposition}

\newtheorem*{corollary*}{Corollary}

\newtheorem*{result*}{Result}

\theoremstyle{definition}
\newtheorem*{definition*}{Definition}

\newtheorem*{example*}{Example}

\newtheorem*{remark*}{Remark}
\newtheorem*{remarks*}{Remarks}

\numberwithin{equation}{section}

\def\rh{\rho}

\def\ph{\phi}

\def\Si{\Sigma}

\def\N{\mathbb{N}}

\def\R{\mathbb{R}}

\def\cG{\mathcal{G}}

\def\cO{\mathcal{O}}

\def\cU{\mathcal{U}}
\def\cV{\mathcal{V}}

\def\sr#1%
{\ifmmode{}^\dagger\else${}^\dagger$\fi\ifvmode
\vbox to 0pt{\vss
 \hbox to 0pt{\hskip\hsize\hskip1em
 \vbox{\hsize3cm\eightpoint\raggedright\pretolerance10000
 \noindent #1\hfill}\hss}\vss}\else
 \vadjust{\vbox to0pt{\vss%
 \hbox to 0pt{\hskip\hsize\hskip1em%
 \vbox{\hsize3cm\eightpoint\raggedright\pretolerance10000%
 \noindent #1\hfill}\hss}\vss}}\fi%
}

\renewcommand{\o}{\circ}

\def\<{\langle}
\def\>{\rangle}

\def\rr{\rightrightarrows}
\def\-{\backslash}

\title{Orbit projections of proper Lie groupoids as fibrations}
 
\author{Armin Rainer}

\address{Armin Rainer: Fakult\"at f\"ur Mathematik, Universit\"at Wien,
Nordbergstrasse~15, A-1090 Wien, Austria}

\email{armin.rainer@univie.ac.at}

\begin{document}

\thanks{The author was supported by 
`Fonds zur F\"orderung der wissenschaftlichen Forschung, Projekt P 17108 N04 \& Projekt P19392'}
\date{December 3, 2007}

\maketitle

\vspace{-0.5cm}

\begin{abstract}
Let $\cG \rr M$ be a source locally trivial proper Lie groupoid such that each orbit is of finite type.
The orbit projection $M \to M/\cG$ is a fibration if and only if $\cG \rr M$ is regular.
\end{abstract}

\vspace{0.5cm}

\cite[2.3]{orb-fib} states that the orbit projection $M \to M/G$ of a proper $G$-manifold $M$ is a \emph{fibration},
i.e., has the homotopy lifting property, 
if and only if $M$ is regular, i.e., each connected component of $M$ has only one orbit type. 
We generalize this to proper Lie groupoids:

\begin{theorem*}
Let $\cG \rr M$ be a source locally trivial proper Lie groupoid such that each orbit is of finite type.
The orbit projection $M \to M/\cG$ is a fibration if and only if $\cG \rr M$ is regular.
\end{theorem*}

Consider a Lie groupoid $\cG \rr M$ with source map $s$ and target map $t$. 
All manifolds in this note are assumed to be smooth, Hausdorff, paracompact, and finite dimensional.
There exist interesting non-Hausdorff groupoids $\cG$, however,
unlike convention we assume that {\it all our groupoids $\cG$ are Hausdorff}.
The Lie groupoid $\cG \rr M$ is said to be \emph{source locally trivial} if the source map $s : \cG \to M$ makes $\cG$ into a 
locally trivial fibration.
For $x \in M$ the set $\cO_x := \{t(g) : g \in s^{-1}(x)\}$ is called \emph{orbit} through $x$.
Each orbit $\cO$ is a regular immersed submanifold in $M$. It is called of \emph{finite type} if there is a 
proper function $f : \cO \to \R$ with a finite number of critical points.  
The connected components of the orbits constitute the leaves of a singular foliation of $M$.
The quotient $M/\cG := \{\cO_x : x \in M\}$ is called the \emph{orbit space} of $\cG \rr M$.
The \emph{orbit projection} $M \to M/\cG$ is the canonical mapping $x \mapsto \cO_x$.
A Lie groupoid $\cG \rr M$ is called \emph{proper} if $(s,t) : \cG \to M \times M$ is a proper mapping. 
For proper Lie groupoids all orbits are closed submanifolds, the orbit space is Hausdorff and paracompact, and all 
isotropy groups $\cG_x := s^{-1}(x) \cap t^{-1}(x)$ are compact Lie groups.  
For any subset $L \subseteq M$ we denote by $\cG_L := \{g \in \cG : s(g),t(g) \in L\}$ the restriction of $\cG$ to $L$.

A typical example of a proper Lie groupoid is the action groupoid $G \ltimes M$ associated to a 
proper $G$-manifold $M$: It is the groupoid $G \times M \rr M$ with $s(g,x)=x$ and $t(g,x)=g.x$, 
where $G \times M \to M \times M, (g,x) \mapsto (x,g.x)$ is a proper mapping. 
The orbits and isotropy groups of $G \ltimes M$ coincide with the usual orbits and isotropy groups of the 
action $G \times M \to M$.

The statement of \cite[2.3]{orb-fib} mentioned earlier is the content of the above theorem when restricted to proper 
action groupoids $G \ltimes M$.
The main ingredient in the proof of \cite[2.3]{orb-fib} is the slice theorem for proper actions due to Palais 
\cite{palais}. Weinstein \cite{weinstein} and Zung \cite{zung} proved the following slice theorem for 
source locally trivial proper Lie groupoids $\cG \rr M$: Let $\cO$ be an orbit of finite type. Then there is an
invariant neighborhood $\cU$ of $\cO$ in $M$ such that the restriction 
$\cG_{\cU}$ of $\cG$ to $\cU$ is isomorphic to the restriction of 
$\cG_{\cO} \ltimes N_{\cO}$ to a tubular neighborhood of the zero section in $N \cO$. 
Here $\cG_{\cO}$ is a transitive Lie groupoid over $\cO$ which acts linearly 
on the normal vector bundle $N \cO := T M/T \cO$ of $\cO$ in $M$.

Let $\Si$ be a \emph{slice} at $x \in \cO$, i.e., a submanifold $\Si \subseteq M$ with $\Si \cap \cO =\{x\}$ and 
$T_xM=T_x\Si \oplus T_x\cO$. The restriction $\cG_\Si$ is a proper Lie groupoid with fixed point $x$ and, thus, 
it is locally isomorphic to the linear action groupoid $\cG_x \ltimes T_x \Si$ (see \cite[2.3]{zung}).
We call a source locally trivial proper Lie groupoid $\cG \rr M$ \emph{regular} if, for every $x \in M$, the action 
$\cG_x \times T_x \Si \to T_x \Si$ is trivial. 

\proof
Since $M/\cG$ is paracompact, 
the projection $\pi : M \to M/\cG$ is a fibration if and only if it is a 
local fibration, i.e., each point in $M/\cG$ has a neighborhood $U$ such that 
$\pi|_{\pi^{-1}(U)} : \pi^{-1}(U) \to U$ is a fibration (e.g.\ \cite[XX 3.6]{dugundji}).

Assume that $\pi: M \to M/\cG$ is a fibration. Let $x \in M$, $\cO$ the orbit through $x$, and $\Si$ a slice at $x$.
The restriction $\cG_\Si$ is isomorphic to the restriction of the action groupoid $\cG_x \ltimes T_x \Si$ to a 
neighborhood of zero, by \cite[2.3]{zung}. 
By source local triviality, the neighborhood can be chosen to be invariant (see \cite[3.3]{weinstein}).
Since $t$ is a submersion, $t^{-1}(\Si)$ is a closed submanifold of $\cG$. 
Then the compact Lie group $\cG_x$ acts on $\Si$ and it acts freely on $t^{-1}(\Si)$ (by left translations via the embedding of $\cG_x$ into 
the group of bisections of $\cG_\Si$). 
The submersion $t|_{t^{-1}(\Si)} : t^{-1}(\Si) \to \Si$ is $\cG_x$-equivariant.
By \cite[5.1]{weinstein} (see also \cite[9.1 step 2]{weinstein}), 
there is a $\cG_x$-invariant neighborhood $\cV$ of $t^{-1}(x)$ in $t^{-1}(\Si)$ and a $\cG_x$-equivariant retraction 
$\rh : \cV \to t^{-1}(x)$ such that $(\rh,t|_{\cV}) : \cV \to t^{-1}(x) \times \Si$ is a diffeomorphism
(after possibly shrinking $\Si$). 
Putting $\cU := s(\cV)$, we obtain a retraction $\cU \to \cO$. 
The mappings $s|_{\cV} : \cV \to \cU$ and $s|_{t^{-1}(x)} :  t^{-1}(x) \to \cO$ are principal bundles with 
structure group $\cG_x$ (acting from the left).

\[
\xymatrix{
&&& \cV \ar[d]^{s|_{\cV}} \ar[r] & t^{-1}(x) \ar[d]^{s|_{t^{-1}(x)}} \\
X \times \{0\} \ar[d] \ar[drr]^(0.57){f} \ar[urrr]^{\tilde{f}} &&& \cU \ar[d]^{\pi|_{\cU}} \ar[r] & \cO \\
X \times I \ar[drr]_{\ph} \ar[urrr]^{\bar \ph} \ar[uurrr]^{\tilde{\bar \ph}} 
\ar@{.>}[rr]_(0.57){t \o \tilde{\bar \ph}} && \Si \ar[d]^{\pi|_{\Si}} \ar@{^(->}[ur] & \cU/\cG_{\cU} & \\
&& \cU/\cG_{\cU} \ar@{=}[ur] &&
}
\]

By assumption, $\pi|_{\cU} : \cU \to \cU/\cG_{\cU}$ is a fibration.
We claim that also $\pi|_{\Si} : \Si \to \cU/\cG_{\cU} \cong \Si/\cG_{\Si}$ is a fibration. 
Let $f : X \times \{0\} \to \Si$ be continuous and let $\ph : X \times I \to \cU/\cG_{\cU}$ be a homotopy of 
$\pi|_{\Si} \o f$. Since $\pi|_{\cU} : \cU \to \cU/\cG_{\cU}$ is a fibration, there exists a 
homotopy $\bar \ph : X \times I \to \cU$ of $f$ covering $\ph$. 
We may lift $f$ to a mapping $\tilde f : X \times \{0\} \to \cV$, 
by setting $\tilde f := u \o f$, 
with $u : M \to \cG$ the unit map which sends $z$ to the identity arrow $1_z$ at $z$.
Since $s|_{\cV} : \cV \to \cU$ is a fibration, there exists a homotopy 
$\tilde{\bar \ph} : X \times I \to \cV$ of $\tilde f$ covering $\bar \ph$. 
It follows that $t \o \tilde{\bar \ph} : X \times I \to \Si$ 
is a homotopy of $f$ covering $\ph$. Hence the claim is proved.

But $\cG_\Si$ is isomorphic to the restriction of the action groupoid $\cG_x \ltimes T_x \Si$ (which is just a group action) to a neighborhood of zero, 
and so the fact that its 
orbit projection $\pi|_{\Si}$ is a fibration implies that $x$ is a regular point (by \cite[2.3]{orb-fib}).

Suppose that $\cG \rr M$ is regular. Let $x \in M$ and use the above terminology. 
Then the projection $\pi|_{\cU} : \cU \to \cU/\cG_{\cU}$ identifies with $\cO \times \Si \to \Si$, which obviously 
is a fibration. So $\pi : M \to M/\cG$ is a local fibration and, thus, a fibration.
\endproof

\begin{remarks*}
(1) The assumption of source local triviality cannot be omitted without substitution, as shown by the following example 
(from \cite[3.4]{weinstein}):
Let $M=\R^2 \- \{0\}$ and let $\cG \rr M$ be the equivalence relation on $M$ consisting of all pairs of points lying on the same vertical line. It is proper and all isotropy groups are trivial, but the source map is not locally trivial over any point on the vertical line through $0$. The orbit projection $M \to M/\cG \cong \R$ is no fibration, since the fiber over $0$ has a different homotopy type than the other fibers.

(2) Nor can the assumption that the orbits are of finite type be omitted (see \cite[5.2]{weinstein}):
Let $M=\{(x,n) \in (-1,1) \times \N : n x^2 < 1\}$ and let $f: M \to (-1,1)$ be the projection on the first factor. 
The fiber product $\cG=M \times_{(-1,1)} M$ is a proper Lie subgroupoid of $M \times M$. 
All isotropy groups are trivial. Its orbit projection, 
which identifies with $f : M \to (-1,1)$, is no fibration, since $f^{-1}(0) \cong \N$ while the other fibers are finite.
\end{remarks*}


\begin{thebibliography}{99}


\bibitem{dugundji}
J. Dugundji,
\emph{Topology},
Allyn and Bacon, Inc., Boston, 1966.  
Zbl 0144.21501

\bibitem{palais}
R.S. Palais, 
\emph{On the existence of slices for actions of non-compact Lie groups}, 
Ann. of Math. (2) \textbf{73} (1961), 295--323. 
Zbl 0103.01802

\bibitem{orb-fib}
A. Rainer,
\emph{Orbit projections as fibrations},
to appear in Czechoslovak Math. J., arXiv: math.DG/0610513.

\bibitem{weinstein}
A. Weinstein,
\emph{Linearization of regular proper groupoids},
J. Inst. Math. Jussieu \textbf{1} (3) (2002), 493--511. 
Zbl 1043.58009

\bibitem{zung}
N.T. Zung, 
\emph{Proper groupoids and momentum maps: linearization, affinity, and convexity},
Ann. Sci. Éc. Norm. Supér. (4) \textbf{39} (2006), no. 5, 841--869. 
Zbl pre05137698

\end{thebibliography}
\end{document}